\documentclass[bst/sn-mathphys]{sn-jnl}

\usepackage{calc}
\usepackage{cleveref}
\usepackage{enumitem}
\usepackage{physics}
\usepackage{algpseudocode}
\usepackage{bm}
\usepackage{mathtools}
\usepackage{multirow,multicol}
\usepackage{array}

\Crefname{equation}{}{}
\Crefname{figure}{Fig.}{Figs.}
\Crefname{section}{Sect.}{Sects.}

\jyear{2022}%

\theoremstyle{thmstyleone}%
\theoremstyle{thmstyletwo}%
\newtheorem{remark}{Remark}%

\theoremstyle{thmstylethree}%

\newcommand{\otherOS}{\Psi}
\newcommand{\dtfact}{\delta} 
\newcommand{\timeratio}{\Gamma} 
\newcommand{\Tstrang}{\tau^{S}} 
\newcommand{\Tother}{\tau^\otherOS} 
\newcommand{\tstrang}{\tilde{\tau}^{S}} 
\newcommand{\tother}{\tilde{\tau}^\otherOS} 
\newcommand{\yy}{\mathbf{y}}
\newcommand{\yyy[1]}{\mathbf{\mathbf{y}}^{[#1]}}
\newcommand{\Dt}{\Delta{t}}
\newcommand{\F}{\mathcal{F}}
\newcommand{\Fl[1]}{\mathbf{\mathcal{F}}^{[#1]}}
\newcommand{\Nop}{N}
\newcommand{\pphi}[2]{\varphi^{[#1]}_{#2}}
\newcommand{\bmalpha}{\bm{\alpha}}
\newcommand{\aalpha}[2]{\alpha^{[#1]}_{#2}}
\def\aak1{AK 3-2(i)}
\def\aaak2{AK 3-2(ii)}
\def\ak52{AK 5-2}
\def\os3237{OS 32-3(7)}
\def\aaaak11{AK 11-4}

\newcommand{\tilt}{\tilde{t}}

\newcommand{\rev}{\textcolor{black}}

\begin{document}

\title[Practical assessment of 3-splitting]{\rev{Beyond Strang: A
    practical assessment of some second-order 3-splitting methods}}


	\author*[1,2]{\fnm{Raymond J.} \sur{Spiteri}}\email{spiteri@cs.usask.ca}

	\author[1,3]{\fnm{Arash} \sur{Tavassoli}}\email{art562@usask.ca}
	\equalcont{These authors contributed equally to this work.}

	\author[1,2]{\fnm{Siqi} \sur{Wei}}\email{siqi.wei@usask.ca}
	\equalcont{These authors contributed equally to this work.}

    \author[3]{\fnm{Andrei}
      \sur{Smolyakov}}\email{andrei.smolyakov@usask.ca}

	\affil*[1]{\orgdiv{Department of Computer Science},
      \orgname{University of Saskatchewan}, \orgaddress{\street{110
          Science Place}, \city{Saskatoon}, \postcode{S7N 5K9},
        \state{SK}, \country{Canada}}}

	\affil[2]{\orgdiv{Department of Mathematics and Statistics},
      \orgname{University of Saskatchewan}, \orgaddress{\street{106
          Wiggins Road},
        \city{Saskatoon}, \postcode{S7N 5E6}, \state{SK},
        \country{Canada}}}

	\affil[3]{\orgdiv{Department of Physics and Engineering Physics},
      \orgname{University of Saskatchewan}, \orgaddress{\street{116
          Science Place}, \city{Saskatoon}, \postcode{S7N 5E2},
        \state{SK}, \country{Canada}}}

    \abstract{Operator splitting is a popular divide-and-conquer
      strategy for solving differential equations. Typically, the
      right-hand side of the differential equation is split into a
      number of parts that are then integrated separately. Many
      methods are known that split the right-hand side into two
      parts. This approach is limiting, however, and there are
      situations when 3-splitting is more natural and ultimately more
      advantageous. The second-order Strang operator-splitting method
      readily generalizes to a right-hand side splitting into any
      number of operators. It is arguably the most popular method for
      3-splitting because of its efficiency, ease of implementation,
      and intuitive nature.  Other 3-splitting methods exist, but they
      are less well-known, and \rev{analysis and} evaluation of their
      performance in practice are scarce. We demonstrate the
      effectiveness of some alternative 3-split, second-order methods
      to Strang splitting on two problems: the reaction-diffusion
      Brusselator, which can be split into three parts that each have
      closed-form solutions, and the kinetic Vlasov--Poisson equations
      that is used in semi-Lagrangian plasma simulations. We find
      alternative second-order 3-operator-splitting methods that
      realize efficiency gains of 10\%--20\% over traditional Strang
      splitting.  \rev{Our analysis for the practical assessment of
        efficiency of operator-splitting methods includes the
        computational cost of the integrators and can be used in
        method design.}  }


\maketitle

\section{Introduction}\label{sec:intro}

The mathematical modeling of the evolution of natural systems is
commonly performed by means of differential equations.
These differential equations often have contributions from
distinct physical processes. In such cases, it has proven to be
computationally fruitful to treat the individual terms (or
\textit{operators}) separately, i.e., by splitting methods.  Different
integrators can be used that take advantage of the specific properties
of the individual operators. For example, implicit-explicit (IMEX)
methods split the right-hand side into two operators, ideally one
stiff and one non-stiff, and treat them with an implicit and an
explicit method, respectively~\cite{ars1997}. When the right-hand side
consists of contributions from many processes, however, IMEX methods
can suffer from the challenge of how to split the right-hand side into
only two operators. In this study, we consider two systems, the
Brusselator system~\cite{Lefever1971} and the Vlasov--Poisson equations,
that admit useful 3-splittings.

Because they are less frequently used, 3-splitting methods are less
well known than 2-splitting methods. The well-known second-order
Strang--Marchuk (or simply \textit{Strang}) splitting
method~\cite{Strang1968, Marchuk1971} admits an intuitive
generalization to any number of operators, and because of this, it is
arguably the most popular 3-splitting method. Other 3-splitting
methods can be derived by considering the order conditions of the
splitting method together with favorable properties such as symmetry
and local error minimization \cite{auzinger2016practical,
  os_coeff_web}. Moreover, adaptive 3-splitting methods and
3-splitting methods with complex coefficients have been proposed
\cite{os_coeff_web}. These 3-splitting methods are less well-known,
less intuitive, and their performance in practice has not been widely
studied.

In \cite{auzinger2017}, the authors compared an adaptive third-order
2-splitting method PP\,3/4\,A\,c~\cite{os_coeff_web} with an adaptive
third-order 3-splitting method PP\,3/4\,A\,3\,c~\cite{os_coeff_web}
applied to the Gray--Scott equations. They concluded that although the
3-splitting method requires fewer steps for a given accuracy, the
additional computation cost required per step ultimately led to
underperformance.

For some problems, however, a 3-splitting method can lead to more
operators that have a closed-form solution that can be evaluated
efficiently. In \cite{crouseilles2015hamiltonian}, the Vlasov--Maxwell
system is split into three parts based on its Hamiltonian.  Although
this splitting allowed the two parts to be solved exactly in time, the
rotation part was solved approximately using Strang splitting. This
approach is revisited in \cite{bernier2020}, where the authors split
the Vlasov--Maxwell equations into three parts again, but instead of
approximating the rotation operator with a Strang splitting method,
they developed a new splitting method to solve the rotation part
exactly in time.  Therefore, the three parts of the Vlasov--Maxwell
system are solved exactly in time, significantly reducing the error of
the overall solution and decreasing the computational cost compared to
using a 2D interpolation.  They also applied several 3-splitting
methods of orders $2$, $3$, $4$, and $6$ to the Vlasov--Maxwell
system, including \aaak2 and \ak52 from~\cite{os_coeff_web}. For
3-splitting methods of order two, it was observed that \aaak2
generates the smallest error and hence is the most efficient
second-order method.

The goal of this paper is to examine the performance in terms of
overall computation efficiency of alternatives to Strang splitting for
problems that lend themselves well to 3-splitting. A key observation
is that although some second-order methods may yield smaller errors
for a given step size, their additional computational cost per step
may negate their overall efficiency.

The outline of the remainder of the paper is as follows. Some general
theoretical background on operator splitting methods along with the
specific operator-splitting methods considered in this study are given
in \cref{sec:theory}. Descriptions of the problems used to illustrate
the performance of the methods follow in~\cref{sec:probs}, and the
performance results themselves appear in~\cref{sec:results}. Finally,
some discussion and conclusions of the study are given
in~\cref{sec:conclusions}.

\section{Theoretical Background}\label{sec:theory}

In this section, we present some background on operator-splitting
methods as discussed in \cite{hairer2006}, with a focus on
$3$-splitting methods.

Consider the initial-value problem (IVP) for an $\Nop$-additively
split ordinary differential equation
\begin{equation} \label{cauchy_problem}
	\dv{\yy}{t} = \F(t,\yy) = \sum\limits_{\ell=1}^\Nop \Fl[\ell](t,\yy),
	\qquad  \yy(0) = \yy_0.
\end{equation}
Let $\pphi{\ell}{\Dt}$ be the flow of the sub-system
\begin{equation} \label{eq:subsys}
\dv{\yyy[\ell]}{t} = \Fl[\ell](t, \yyy[\ell])
\end{equation}
for $\ell=1,2,\dots,\Nop$. \rev{We refer to the $\pphi{\ell}{\Dt}$
  generally as \emph{sub-integrations} because they may be
  approximated numerically.} Compositions of $\pphi{\ell}{\Dt}$ can be
used to construct numerical solutions to \cref{cauchy_problem}. For
example, the following two methods
\begin{subequations}
	\begin{align}
		& \Phi_{\Dt} :=  \pphi{\Nop}{\Dt} \circ \pphi{\Nop-1}{\Dt}\circ \cdots \circ \pphi{1}{\Dt}, \label{eq:os_godunov} \\
		& \Phi_{\Dt}^\ast :=  \pphi{1}{\Dt} \circ \pphi{2}{\Dt}\circ \cdots \circ \pphi{\Nop}{\Dt}  \label{eq:os_godunovadj}
	\end{align}
\end{subequations}
are commonly known as the Godunov (or Lie--Trotter) splitting
methods. The two methods \cref{eq:os_godunov} and
\cref{eq:os_godunovadj} are adjoints of each other and are both
first-order accurate. Moreover, the flows
$\left\{\pphi{\ell}{\Dt}\right\}_{\ell=1}^\Nop$ can be composed in any
order to create a first-order accurate Godunov splitting method
$\Phi_{\Dt}^{G}$, and its adjoint $\left.\Phi_{\Dt}^{G}\right.^{\ast}$
can be derived by reversing the order of composition.

Another popular method is the Strang splitting method, which is
constructed by composing $\Phi_{\Dt}^{G}$ and
$\left.\Phi_{\Dt}^{G}\right.^{\ast}$ with halved step sizes. Using
$\cref{eq:os_godunov}$ and $\cref{eq:os_godunovadj}$, we can write a
Strang splitting method as
\begin{equation} \label{strang}
\begin{aligned}
	\Psi_{\Dt}^{S} & = \left.\Phi_{\Dt/2}^{G}\right.^{ \ast} \circ \Phi^G_{\Dt/2} \\
	& = \pphi{1}{\Dt/2} \circ \pphi{2}{\Dt/2}\circ \cdots
	\circ  \pphi{\Nop-1}{\Dt/2}\circ  \pphi{\Nop}{\Dt} \circ \pphi{\Nop-1}{\Dt/2}\circ \cdots \circ \pphi{1}{\Dt/2}.
\end{aligned}
\end{equation}
We note that the terms $ \pphi{\Nop}{\Dt/2}\circ \pphi{\Nop}{\Dt/2} $
can be combined together as $ \pphi{\Nop}{\Dt}$ by group property if
$\pphi{\Nop}{\Dt/2}$ is the exact flow.  If $\pphi{\Nop}{\Dt/2}$ is
approximated, however, these two methods are different numerical
methods with different accuracy and stability properties.

We express the general form of an $s$-stage operator-splitting methods
as follows. Let
$\bmalpha = [\bmalpha_1, \bmalpha_2, \dots, \bmalpha_s]$, where
$\bmalpha_k = [\aalpha{1}{k}, \aalpha{2}{k}, \dots ,
\aalpha{\Nop}{k}], \ k=1,2,\dots, s$, be coefficients of an
operator-splitting method. An $s$-stage operator-splitting method that
solves \eqref{cauchy_problem} can be written as
\begin{equation} \label{eq:os_method}
	\Psi_{\Dt} := \prod_{k=1}^{s} \Phi_{\bmalpha_k \Dt}^{\{k\}} = \Phi_{\bmalpha_s \Dt}^{\{s\}} \circ \Phi_{\bmalpha_{s-1} \Dt}^{\{s-1\}} \circ \cdots \circ \Phi_{\bmalpha_1 \Dt}^{\{1\}},
\end{equation}
where
$\Phi_{\bmalpha_k \Dt}^{\{k\}} := \pphi{\Nop}{\aalpha{\Nop}{k}\Dt}
\circ\pphi{\Nop-1}{\aalpha{\Nop-1}{k}\Dt} \circ \cdots \circ
\pphi{1}{\aalpha{1}{k}\Dt} $.

Henceforth, we focus on $3$-splitting methods; i.e., $\Nop =3$. For
the numerical method \cref{eq:os_method} to have order $p_{OS}$, the
operator-splitting coefficients $\bmalpha$ must satisfy a system of
order conditions in the form of polynomial equations. These equations
can be derived from the well-known Baker–Campbell–Hausdorff (BCH)
formula; see, e.g.,~\cite{hairer2006}. The approach for deriving
operator-splitting methods described in \cite{Auzinger2014} also
relies on the BCH formula, but the order conditions are generated
automatically via computer algebra. Because our focus is on examining
alternatives to Strang splitting with 3 operators, we focus on order
conditions up to $p_{OS}=2$ for the $3$-splitting case:

\begin{subequations}
\tiny{
  \begin{align}
    \label{order1}
    p_{OS} &= 1: &&\quad \sum\limits_{k=1}^s \aalpha{1}{k} = 1,
    &&\quad \sum\limits_{k=1}^s \aalpha{2}{k} = 1,
    && \quad \sum\limits_{k=1}^s \aalpha{3}{k} = 1, \\
    \label{order2}
    p_{OS} &= 2:  &&\quad \sum\limits_{k=1}^s \aalpha{1}{k}\left(\sum\limits_{k'=k}^s \aalpha{2}{k'}\right)  = \frac{1}{2},
    && \quad \sum\limits_{k=1}^s \aalpha{1}{k}\left(\sum\limits_{k'=k}^s \aalpha{3}{k'}\right)  = \frac{1}{2},
    && \quad \sum\limits_{k=1}^s \aalpha{2}{k}\left(\sum\limits_{k'=k}^s \aalpha{3}{k'}\right)  = \frac{1}{2}.
  \end{align}}

\end{subequations}

In this study, we are specifically interested in applications where
the sub-systems \cref{eq:subsys} can be solved exactly or with high
precision. That is, the main source of error in the numerical solution
is the splitting error. In \cite[\rev{equation
  (4.2b)}]{auzinger2016practical}, the authors developed a local error
measure (LEM) based on the operator-splitting coefficients
$\bmalpha$. Four second-order 3-splitting methods are proposed in
\cite{os_coeff_web}. They are Strang, \aak1, \aaak2, and \ak52. The
efficiency of an operator-splitting method is affected by both the
number of sub-integrations required by the method and (to some extent)
by the LEM. To compare with Strang, we consider the class of
three-stage, second-order 3-splitting methods denoted by OS 32-3. To
satisfy the order conditions \cref{order1} and \cref{order2}, an OS
32-3 method requires at least 5 sub-integrations. \rev{For an OS 32-3
  method to have exactly 5 sub-integrations, 4 of the 9 coefficients
  $\{\aalpha{\ell}{k}\}_{k,l = 1,2,3}$ must be zero. Order condition
  \cref{order1} indicates that, for each operator $\ell$, at most two
  of $\{\aalpha{\ell}{k}\}_{k=1,2,3}$ can be zero. For each
  $\ell = 1,2,3$, setting two of $\{\aalpha{\ell}{k}\}_{k=1,2,3}$ and
  one of $\{\aalpha{\ell'}{k}\}_{k=1,2,3, \ell' \neq \ell}$ equal to
  zero, we can solve the system \cref{order1} and \cref{order2} to
  show that any OS 32-3 method that requires exactly 5
  sub-integrations is equivalent to Strang splitting with different
  permutations of the operators $\ell$. For example,
  $\pphi{1}{\Dt/2} \circ \pphi{2}{\Dt/2}\circ \pphi{3}{\Dt} \circ
  \pphi{2}{\Dt/2}\circ \pphi{1}{\Dt/2}$,
  $\pphi{2}{\Dt/2} \circ \pphi{3}{\Dt/2}\circ \pphi{1}{\Dt} \circ
  \pphi{3}{\Dt/2}\circ \pphi{2}{\Dt/2}$, and
  $\pphi{3}{\Dt/2} \circ \pphi{1}{\Dt/2}\circ \pphi{2}{\Dt} \circ
  \pphi{1}{\Dt/2}\circ \pphi{3}{\Dt/2}$ are three of the six possible
  OS 32-3 Strang-splitting methods with exactly 5 sub-integrations.  }

\aak1 is the OS 32-3 method with the smallest LEM for 6
sub-integrations. We constructed an OS 32-3 method with the smallest
LEM that requires 7 sub-integrations, \rev{but it was not competitive
  and is not considered further}. \Cref{tab:os3_schemes} summarizes
the main characteristics of these five methods. The operator-splitting
coefficients of \aak1, \aaak2, and \ak52 are given in
\cref{tab:ak32i}, \cref{tab:ak32ii}, and \cref{tab:ak52},
respectively.
\begin{table}[htbp]
	\centering
	\begin{tabular}{|c|c|c|c|}
      \hline
      Method &  Stages & Sub-integrations & LEM \\
      \hline
      Strang & 3 & 5 & 1.48  \\ 	\hline
      \aak1 & 3 & 6 & 1.06 \\	\hline
      \aaak2 & 3 & 9 & 0.29 \\	\hline
      \ak52 & 5 & 9 & 0.22 \\
      \hline
	\end{tabular}
	\caption{Summary of second-order 3-splitting methods considered.
      \label{tab:os3_schemes}}
\end{table}

\begin{table}[htbp]
	\centering
	\begin{tabular}{|c|c|c|c|}
		\hline
		$k$&  $\aalpha{k}{1}$ & $\aalpha{k}{2}$ & $\aalpha{k}{3}$ \\
		\hline
		1 & $0.5$ & $1-1/\sqrt{2}$ & $1/\sqrt{2}$\\
		\hline
		2 & $0$ & $1/\sqrt{2}$ & $1-1/\sqrt{2}$\\
		\hline
		3 & $0.5$ & $0$ & $0$\\
		\hline
	\end{tabular}
	\caption{Operator-splitting coefficients of \aak1
		\label{tab:ak32i}}
\end{table}


\begin{table}[htbp]
	\centering
	\begin{tabular}{|c|c|c|c|}
		\hline
		$k$&  $\aalpha{k}{1}$ & $\aalpha{k}{2}$ & $\aalpha{k}{3}$ \\
		\hline
		1 & $0.316620935432115636$ & $0.273890572734778059$ & $0.662265355057626845$\\
		\hline
		2 & $-0.0303736077786568570$ & $0.438287559165397521$ & $0.0664399910533392230$\\
		\hline
		3 & $0.713752672346541221$ & $0.287821868099824420$ & $0.271294653889033932$\\
		\hline
	\end{tabular}
	\caption{Operator-splitting coefficients of \aaak2
		\label{tab:ak32ii}}
\end{table}

\begin{table}
	\centering
	\begin{tabular}{|c|c|c|c|}
		\hline
		$k$&  $\aalpha{k}{1}$ & $\aalpha{k}{2}$ & $\aalpha{k}{3}$ \\
		\hline
		1 & $0.161862914279624$ & $0.242677859055102$ & $0.5$ \\
		\hline
		2 &  $0.338137085720376$ & $0.514644281889796$ & $0$ \\
		\hline
		3 & $0.338137085720376$ & $	0 $ & $	0.5$ \\
		\hline
		4 & $0$ & $0.242677859055102$ & $0$ \\
		\hline
		5 & $0.161862914279624$ & $	0$ & $0$ \\
		\hline
	\end{tabular}
	\caption{Operator-splitting coefficients of \ak52
		\label{tab:ak52}}
\end{table}

\rev{\subsection{Practical assessment of splitting methods}}

\rev{To assess the overall efficiency gain of a splitting method
$\otherOS$ against the Strang splitting method, we must consider not
only the total number of time steps required by each method but also
the cost per step. Let $N^{S}$ and $\tstrang$ ($N^{\otherOS}$ and
$\tother$) be the numbers of time steps and the wall clock times per
time step of the Strang ($\otherOS$) splitting method. The efficiency
gain $\eta$ of a splitting method $\otherOS$ with respect to Strang
splitting is defined in terms of work-precision as
\begin{equation*}
	\eta \coloneqq \frac{\Tstrang-\Tother}{\Tstrang},
\end{equation*}
where $\Tstrang$ and $\Tother$ are the total wall clock time
taken by the Strang and $\Psi$ splitting methods, respectively.
Hence,
\begin{equation}\label{eq:efficiency1}
	\eta = \frac{N^S\tstrang - N^\otherOS \tother}{N^S\tstrang}
	= 1 - \frac{N^\otherOS}{N^S} \frac{\tother}{\tstrang}
 = 1 - \frac{\Dt^S}{\Dt^\otherOS} \frac{\tother}{\tstrang},
\end{equation}
where $\Dt^S$ and $\Dt^\otherOS$ are the largest step sizes for Strang
and $\otherOS$ splitting method to achieve certain accuracy,
respectively. The ratio
$\displaystyle \displaystyle \frac{\tother}{\tstrang}$ can be
estimated using the cost of additional sub-integrations required by
performing $\otherOS$ compared with Strang. As indicated in
\cref{sec:theory}, the Strang splitting method requires $k^S = 5$
sub-integrations. Let $k^\otherOS$ be the number of sub-integrations
required by the $\otherOS$ splitting method, then
\begin{equation}\label{eq:extracost}
\begin{aligned}
  \frac{\tother}{\tstrang} & \approx \frac{\tstrang + \sum\limits_{j=1}^{k^\otherOS -k^S} \tother_j}{\tstrang}, \\
                           & = 1 + \sum\limits_{j=1}^{k^\otherOS -k^S}\frac{\tother_j}{\tstrang},
\end{aligned}
\end{equation}
where $\tother_j$ is the wall clock time of additional sub-integration
$j$. We note that with this formula, we have taken into account that
the cost of solving each operator can be different. Now we define the
ratio $\displaystyle \dtfact := \frac{\Dt^\otherOS}{\Dt^S}$ and the
extra-time fraction
$\displaystyle \timeratio = \sum\limits_{j=1}^{k^\otherOS
  -k^S}\frac{\tother_j}{\tstrang}$. Using \cref{eq:efficiency1} and
$\cref{eq:extracost}$, we can write $\eta$ as
\begin{equation}\label{eq:efficiency}
	\eta = 1- \frac{1+\timeratio}{\dtfact}.
\end{equation}
}

\begin{remark}
  \rev{Many} benchmarking studies consider the efficiency of a
  numerical method \rev{solely} as its accuracy per time
  step~\cite{bernier2020,casas2020composition}. However, such
  approaches do not consider the computational cost of each step
  required to achieve a given accuracy. In other words, the accuracy
  per time step of two competing methods can be used to calculate
  \rev{the ratio of largest acceptable step sizes} $\delta$, but it
  neglects the \rev{extra-time fraction} $\timeratio$. Similarly, even
  if it were a perfect error estimator, the LEM does not include any
  information about $\displaystyle \timeratio$. Therefore, it
  \rev{should not be} surprising that a particular method can be less
  efficient than another despite having a smaller LEM.
\end{remark}

\section{Problem Descriptions}\label{sec:probs}

In this section, we describe two problems used to illustrate the
performance of the second-order 3-splitting methods described
in~\cref{sec:probs}.

\subsection{Brusselator}\label{subsec:Brusselator_prob}

\rev{The first problem considered is the well-known Brusselator
  equation~\cite{Lefever1971}, which is a set of reaction-diffusion
  PDEs that describes an autocatalytic reaction between two chemical
  species with different rates of diffusion. It is a commonly used
  benchmark problem for numerical methods because it exhibits
  interesting dynamics such as periodic solutions and bifurcations
  that are well understood analytically.} The Brusselator problem can
be solved using 2-splitting method according to reaction and
diffusion, e.g.,~\cite{ropp2005}. When split into two operators, the
reaction operator is non-linear and does not have a closed-form
solution. In this study, we split the Brusselator problem into three
operators. This splitting allows for a closed-form solution to each
operator.

The Brusselator problem is defined as
\begin{subequations} \label{eq:brusselator}
	\begin{align}
		\pdv{T}{t} &= D_1 \pdv[2]{T}{x} + \alpha - (\beta+1)T + T^2C, \\
		\pdv{C}{t} &= D_2 \pdv[2]{C}{x} + \beta T -T^2C,
	\end{align}
\end{subequations}
where $T=T(x,t)$ and $C=C(x,t)$ represent concentrations of different
chemical species. The parameter values are $\alpha = 0.6$, $\beta =2$,
and $\displaystyle D_1= D_2 = \frac{1}{40}$, with boundary conditions
$T(0,t) = T(1,t) = \alpha$ and
$\displaystyle C(0,t) = C(1,t) = \frac{\beta}{\alpha}$ and initial
conditions $T(x,0) = \alpha + x(1-x)$ and
$\displaystyle C(x,0) = \frac{\beta}{\alpha} + x^2(1-x)$.
Equation \Cref{eq:brusselator} is split into three parts as
\begin{subequations}
  \label{eq:brusselator_split}
	\begin{align}
		\label{brusselator_1}
		& \left\{\begin{array}{l}
			\displaystyle \pdv{T^{[1]}}{t} = D_1 \pdv[2]{T^{[1]}}{x}, \\[1ex]
			\displaystyle \pdv{C^{[1]}}{t} = D_2 \pdv[2]{C ^{[1]}}{x},
		\end{array} \right.\\
		\label{brusselator_2}
		& \left\{\begin{array}{l}
			\displaystyle \pdv{T^{[2]}}{t} = \alpha - (\beta+1)T^{[2]}, \\[1ex]
			\displaystyle \pdv{C^{[2]}}{t} = \beta T^{[2]},
		\end{array} \right.\\
		\label{brusselator_3}
		& \left\{\begin{array}{l}
			\displaystyle \pdv{T^{[3]}}{t} = (T^{[3]})^2C^{[3]},  \\[1ex]
			\displaystyle \pdv{C^{[3]}}{t} = -(T^{[3]})^2C^{[3]}.
		\end{array} \right.
	\end{align}
\end{subequations}
Equation \cref{brusselator_1} is solved by first discretizing the PDE
using a second-order finite-difference method. Then
\cref{brusselator_1} is converted to
$$\left\{
\begin{aligned}
	\pdv{T}{t} &=  D_1 M_1 T,  \\
	\pdv{C}{t} &=  D_2 M_2 C,
\end{aligned}
\right. $$ where $M_1$ and $M_2$ correspond to the finite difference
stencil for $\displaystyle \pdv[2]{T}{x}$ and
$\displaystyle \pdv[2]{C}{x}$ respectively.  To solve
\cref{brusselator_3}, we note that
$\displaystyle \pdv{T^{[3]}}{t} + \pdv{C^{[3]}}{t} = 0$. Hence,
$T^{[3]}+C^{[3]} = k$ for some constant $k$ with respect to
$t$. Substituting $C^{[3]} = k-T^{[3]}$ into the first equation of
\cref{brusselator_3}, we can obtain a differential equation in
$T^{[3]}$ only:
\begin{equation}\label{DEofT}
 \pdv{T^{[3]}}{t} = (T^{[3]})^2(k-T^{[3]}).
\end{equation}
Equation \cref{DEofT} can be solved using separation of variables and
partial fraction decomposition.  The exact solution of
\cref{eq:brusselator_split} is then given by

\begin{subequations}\label{brusselator_sol}
	\begin{alignat}{3}
		\label{brusselator_sol1}
		&\left\{
		\begin{aligned}
			T^{[1]}(t+\Dt) &=  T^{[1]}(t)\exp(D_1M_1\Dt), \\
			C^{[1]}(t+\Dt) &=  C^{[1]}(t)\exp(D_2M_2\Dt), \\
		\end{aligned}
		\right.  \\
		\label{brusselator_sol2}
		&\left\{
		\begin{aligned}
			T^{[2]}(t+\Dt) &= \frac{\exp(-\Dt(\beta+1))[T^{[2]}(t) -\alpha + T^{[2]}(t)\beta + \alpha \exp(\Dt(\beta+1))]}{\beta+1},\\
			C^{[2]}(t+\Dt) &=  C^{[2]}(t)+\frac{T^{[2]}(t)\beta + \alpha \beta \Dt}{\beta+1} \\
			& - \frac{\beta \exp(-\Dt(\beta+1))[T^{[2]}(t) - \alpha + T^{[2]}(t) \beta + \alpha \exp(\Dt(\beta+1))]}{(\beta +1)^2},
		\end{aligned}
		\right.   \\
		\label{brusselator_sol3}
		& \left\{
		\begin{aligned}
			& \frac{1}{k^2} \ln \left\vert \frac{T^{[3]}(t+\Dt)}{T^{[3]}(t+\Dt) -k} \right\vert - \frac{1}{kT^{[3]}(t+\Dt)} = \Dt + c_1, \\
			& C^{[3]}(t+\Dt) = k - T^{[3]}(t+\Dt),
		\end{aligned}
		\right.
	\end{alignat}
\end{subequations}
where $\displaystyle c_1 =\frac{1}{k^2}\ln\left\vert \frac{T^{[3]}(t)}{T^{[3]}(t)-k} \right\vert - \frac{1}{kT^{[3]}(t)} $ and $k = T^{[3]}(t) + C^{[3]}(t)$.

We report on the performance using Strang, \aak1, \aaak2, and \ak52 to
solve~\cref{eq:brusselator_split} in
\cref{subsec:Brusselator_results}.

\subsection{Plasma dynamics with kinetic Vlasov--Poisson equations}\label{subsec:ECDI_prob}

\rev{The second problem considered involves simulation of the electron
  cyclotron drift instability (ECDI) from the field of plasma
  physics. The ECDI is a plasma instability that is driven by the
  relative drift velocity of ions and electrons in the presence of a
  magnetic field and has recently received a lot of attention from the
  plasma science
  community~\cite{boeuf2018b,charoy2021interaction,sengupta2020mode,asadi2019numerical,hara2020cross,mandal2020cross,janhunen2018nonlinear,janhunen2018evolution,tavassoli2022nonlinear}. A
  number of factors, including high dimension, high resolution
  requirements, large required simulation domains, and long required
  simulation times, conspire to make this a challenging numerical
  problem.  Simulations that provide meaningful data typically take on
  the order of ten days to run using 32 cores (with 80\% scaling
  efficiency). Efficient time integration methods greatly impact such
  problems, where savings of 10\%--20\% translate to days or weeks of
  CPU time.}


The behavior of the ECDI is governed by the Vlasov--Poisson system of
equations. For one spatial dimension and two velocity dimensions
(1D2V), these equations are
\begin{subequations}\label{ecdi_orgin}
\begin{align}
	\pdv{f_e}{t}+v_{x}\pdv{f_e}{x}-\alpha_1\qty(E_x-v_{z}\alpha_2)\pdv{f_e}{v_{x}}-\alpha_1\qty(\alpha_3+v_{x}\alpha_2)\pdv{f_e}{v_{z}}&=0,\label{eq:evlasov}\\
	\pdv{f_i}{t}+v_{x}\pdv{f_i}{x}+ E_x\pdv{f_i}{v_{x}}&=0,\label{eq:ivlasov}\\
	\pdv{E_x}{x}-\int\limits_{-\infty}^{\infty}\int\limits_{-\infty}^{\infty} f_i\dd{v_{x}}\dd{v_{z}}+\int\limits_{-\infty}^{\infty}\int\limits_{-\infty}^{\infty} f_e\dd{v_{x}}\dd{v_{z}}&=0,\label{eq:ef}
\end{align}
\end{subequations}
where the dependent variables $f_e=f_e(x,v_x,v_z,t)$ and
$f_i=f_i(x,v_x,t)$ are distribution functions of the electrons and ions,
respectively, $E_x=E_x(x,t)$ is the electric field, and
$\alpha_1,\;\alpha_2$, and $\alpha_3$ are parameters that depend on
the ion mass, plasma density, and external fields. The total
energy in this system is defined by
\begin{equation}	\label{eq:total_U}
\resizebox{\textwidth}{!}{
$\begin{aligned}
	U(t)= & \frac{1}{2\alpha_1}\int\limits_0^L \int\limits_{-\infty}^{\infty} \int\limits_{-\infty}^{\infty} (v_x^2+v_z^2)f_e(x,v_x,v_z,t)\, \dd v_x\dd v_z\dd x+\frac{1}{2}\int\limits_0^L \int\limits_{-\infty}^{\infty} v_x^2f_i(x,v_x,t)\, \dd v_x\dd x    \\
	& +\frac{1}{2}\int\limits_0^L  E_x^2\, \dd x +\int\limits_0^t\int\limits_0^L \int\limits_{-\infty}^{\infty}\int\limits_{-\infty}^{\infty} \alpha_3f_e(x,v_x,v_z,\tilt)v_z\, \dd v_x\dd v_z\dd x\, \dd \tilt,
\end{aligned}$}
\end{equation}
where $L$ is the size of the physical domain. Here, we use a periodic
boundary condition for the physical domain. For both $f_e$ and $f_i$,
we have
\begin{equation}
	\lim_{v_x\rightarrow \pm \infty}f_s=0, \
	\lim_{v_z\rightarrow \pm \infty}f_s=0, \quad s=e,i.
\label{eq:bound_vxz}%
\end{equation}
These conditions mean that there can be no particle with an infinite
velocity. Conditions \Cref{eq:bound_vxz} are imposed numerically by
adjusting the boundaries so that $f_s\approx 10^{-8}$ in their
vicinity.

A popular Eulerian--Vlasov numerical method is the
\textit{semi-Lagrangian} method
\cite{cheng1976integration,cheng1977integration,sonnendrucker1999semi,
crouseilles2010conservative,crouseilles2009forward,qiu2010conservative,
besse2008convergence,ghizzo2003non}. In this method,
\Cref{eq:evlasov,eq:ivlasov} are split into a number of sub-equations,
and each sub-equation is integrated using the method of
characteristics. For both the ions and electrons, the characteristic
equations of the Vlasov equation are in the form of equations of
particle motion. In this study, the semi-Lagrangian method, as
discussed below, is used for the integration of the Vlasov--Poisson system.

We first consider the electron Vlasov equation \Cref{eq:evlasov}.
For this equation, the characteristic equations are defined by
\begin{subequations}
  \label{eq:ch}
\begin{align}
	\dv{X_e}{t}=&V_{X_e},\label{eq:ch_xe}\\
	\dv{V_{X_e}}{t}=&-\alpha_1 \qty(E_x(x,t)-V_{Z_e}\alpha_2),\label{eq:ch_vxe}\\
	\dv{V_{Z_e}}{t}=&-\alpha_1 (\alpha_3+V_{X_e}\alpha_2).\label{eq:ch_vze}
\end{align}
\end{subequations}
Now, one can define characteristic curves
$$\vb{W}(t{;}x,v_x,v_z)\equiv
(X_e(t{;}x,v_x,v_z),V_{X_e}(t{;}x,v_x,v_z),V_{Z_e}(t{;}x,v_x,v_z))$$ as
functions that satisfy \Cref{eq:ch} and at time $t+\Delta t$ are terminated
at the point $(x,v_x,v_z)$. From \Cref{eq:evlasov,eq:ch}, one can show
that the distribution function is constant along the characteristic
curves, i.e., $\qty[\dv{f_e}{t}]_{\vb{W}}=0$, and accordingly,
\begin{equation}
	f_e(x,v_x,v_z,t+\Delta t)=f_e(X_e(t{;}x,v_x,v_z),V_{X_e}(t{;}x,v_x,v_z),V_{Z_e}(t{;}x,v_x,v_z),t).
\end{equation}
Therefore, for updating the solution of \Cref{eq:evlasov} at
$t+\Delta t$, one only needs to find the base of the characteristic
curves $\vb{W}(t{;}x,v_x,v_z)$. In general, to find the base of the characteristic
curve, one needs to perform a backward integration. 
Because of this, the $\vb{W}(t{;}x,v_x,v_z)$ is in general found
implicitly by solving a nonlinear equation
$\vb{F}(\vb{W}(t{;}x,v_x,v_z),E_x(t))=\vb{0}$~\cite{coulaud1999parallelization,magdi2009method}. However,
this particular problem can be solved by splitting the Vlasov equation
into equations along the directions of one of the independent
variables at a time. For example, \Cref{eq:evlasov} is split into
three sub-equations
\begin{subequations}
  \label{eq:split}
\begin{align}
	\pdv{f_e^{[1]}}{t}+v_{x}\pdv{f_e^{[1]}}{x}&=0,\label{eq:xe_split}\\
	\pdv{f_e^{[2]}}{t}+a_z\pdv{f_e^{[2]}}{v_{z}}&=0,\label{eq:vze_split}\\
	\pdv{f_e^{[3]}}{t}+a_x\pdv{f_e^{[3]}}{v_{x}}&=0,\label{eq:vxe_split}
\end{align}
\end{subequations}
where $a_z\equiv -\alpha_1\qty(\alpha_3+v_{x}\alpha_2)$ and
$a_x\equiv -\alpha_1\qty(E_x-v_{z}\alpha_2)$ . In \Cref{eq:xe_split},
the coefficient of $\displaystyle \pdv{f_e^{[1]}}{x}$ is $v_x$, which
is a constant of that equation. Similarly in
\Cref{eq:vze_split,eq:vxe_split}, $a_x$ and $a_z$ are
constant. Because of this property of \Cref{eq:split}, the exact base
of the characteristic of these equations can be found explicitly.
Using the method of characteristics, solutions of \Cref{eq:split} at
$t+\Delta t$ are
\begin{subequations}
  \label{eq:shift}
\begin{align}
	f_e^{[1]}(x,t+\Delta t)&= f_e^{[1]}(x-v_x\Delta t),\label{eq:shift_xe}\\
	f_e^{[2]}(v_z,t+\Delta t)&=f_e^{[2]}(v_z-a_z\Delta t),\label{eq:shift_vze}\\
	f_e^{[3]}(v_x,t+\Delta t)&= f_e^{[3]}(v_x-a_x\Delta t)\label{eq:shift_vxe}.
\end{align}
\end{subequations}
We note that \Cref{eq:shift} are the exact solutions of
\Cref{eq:split}. However, because at time $t$, $f_e$ is only known at
a particular grid point $(x_n,v_{xn},v_{zn})$, evaluating
\Cref{eq:shift} requires an interpolation. Many interpolation methods
are suggested in the literature; among them, the cubic spline is one
of the most popular because it is believed to keep a good balance
between the accuracy and the performance in semi-Lagrangian
methods~\cite{staniforth1991semi}. In this study, all the
interpolations are done using cubic splines.

Similar to \Cref{eq:evlasov}, \Cref{eq:ivlasov} is split into two equations
\begin{subequations}
\begin{align}
	\pdv{f_i^{[1]}}{t}+v_{x}\pdv{f_i^{[1]}}{x}&=0,\label{eq:xi_split}\\
	\pdv{f_i^{[2]}}{t}+E_x\pdv{f_i^{[2]}}{v_{x}}&=0\label{eq:vxi_split},
\end{align}
\end{subequations}
which are solved using the method of characteristics
\begin{subequations}
\begin{align}
	f_i^{[1]}(x,t+\Delta t)&= f_i^{[1]}(x-v_x\Delta t),\label{eq:shift_xi}\\
	f_i^{[2]}(v_x,t+\Delta t)&=f_i^{[2]}(v_x-E_x\Delta t).\label{eq:shift_vzi}
\end{align}
\end{subequations}
Although \Cref{eq:ivlasov} is split into two operators, we cast it as
a 3-splitting method with a trivial operator [3], i.e.,~a trivial
operator in the $v_z$ direction. 
Equation \cref{ecdi_orgin} can then split into the following three operators:

\begin{subequations}\label{ecdi_split}
	\begin{align}
		\label{ecdi_1}
		& \left\{\begin{array}{l}
			\displaystyle \pdv{f_e^{[1]}}{t} = -v_x \pdv{f_e^{[1]}}{x}, \\[1ex]
			\displaystyle \pdv{f_i^{[1]}}{t} = -v_x \pdv{f_i^{[1]}}{x},
		\end{array} \right.\\
		\label{ecdi_2}
		& \left\{\begin{array}{l}
			\displaystyle \pdv{f_e^{[2]}}{t} = -a_z \pdv{f_e^{[2]}}{v_z}, \\[1ex]
		\end{array} \right.\\
		\label{ecdi_3}
		& \left\{\begin{array}{l}
			\displaystyle \pdv{f_e^{[3]}}{t} = -a_x \pdv{f_e^{[3]}}{v_x},   \\[1ex]
			\displaystyle \pdv{f_i^{[3]}}{t} = -E_x \pdv{f_i^{[3]}}{v_x} .
		\end{array} \right.
	\end{align}
\end{subequations}

The split Vlasov equations \cref{ecdi_split} are predominantly solved
using the Strang splitting method \cref{strang} for three
operators. Equation \cref{eq:ef} is integrated after the first stage
$\bmalpha_1$ of the operator-splitting process.

%
%

We note that the energy $U(t)$ defined by \Cref{eq:total_U} is
theoretically constant. However, the semi-Lagrangian method is not
energy conserving. Because of this, the so-called phenomena of
``numerical heating" and ``numerical cooling'' are concerns in
semi-Lagrangian simulations. Accordingly, the deviation from energy
conservation is a reasonable indicator of simulation error. This is
especially true when qualitative global aspects of the solution are of
interest more so than specific quantitative local values.

%

 For this problem, we test the \aak1 splitting method in addition
 to the Strang method. Equation \cref{eq:ef} is again integrated after
 the first stage $\bmalpha_1$.  We note that the \aak1 method has one
 more sub-integration compared to the Strang method. Accordingly, it
 generally has a greater computational expense on a per step basis. In
 order for \aak1 to outperform Strang, therefore, it must allow for an
 increased step size that more than offsets this additional cost per
 step while maintaining sufficient accuracy.

\section{Results}\label{sec:results}

In this section, we describe the results from applying a number of
second-order, 3-splitting methods to the Brusselator and ECDI
problems. We find that some alternatives to traditional Strang
3-splitting can be 10\%--20\% more efficient for a given level of
accuracy.

\subsection{Brusselator}\label{subsec:Brusselator_results}

A reference solution for the Brusslator problem~\cref{eq:brusselator}
for $t\in [0,80]$ is computed using the MATLAB parabolic and elliptic
PDE solver {\tt pdepe}. We decreased the spatial meshsize $\Delta x$
and adjusted the absolute and relative tolerances for the solver until
there were at least $6$ matching digits between successive
approximations at 32,000 and 800 uniformly distributed points in space
and time, respectively.  For our experiments with the various
splitting methods, the spatial derivatives are discretized using
central finite differences on a uniform grid on the interval
$x\in [0,1]$.

The numerical experiments using operator-splitting methods are
implemented using Python. Although all sub\rev{-}integrators have
analytical solutions, the solution \cref{brusselator_sol1} requires
the matrix exponential, and the solution \cref{brusselator_sol3}
requires the solution to a non-linear equation. It turns out to be
faster to obtain an ``exact'' solution by solving \cref{brusselator_1}
and \cref{brusselator_3} with the Python {\tt
  scipy.integrate.solve\_ivp} function using {\tt rtol = 1e-10} and
{\tt atol = 1e-13} and using the closed-form for \cref{brusselator_2}.

\begin{figure}[htbp]
	\centering
	\includegraphics[width=5in]{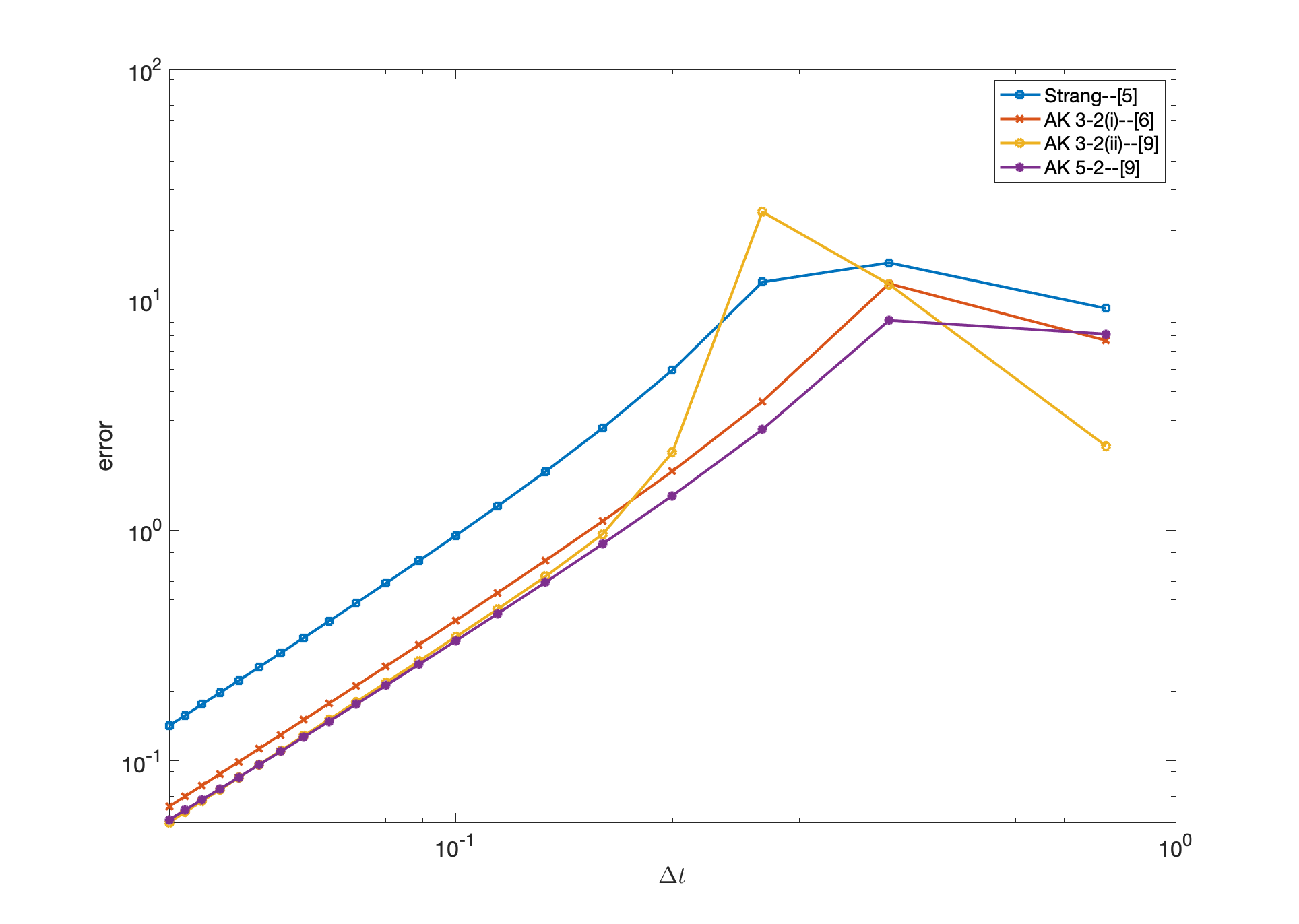}
	\caption{Convergence of second-order splitting methods. All
      methods demonstrate the expected
      convergence rate. \label{fig:brusselator_convergence}}
\end{figure}

In this example, the error in the numerical solution $Y=[T;C]$ is
measured using the mixed root mean square (MRMS) error defined by
\begin{equation*}
  \text{error} = [\text{MRMS}]_Y = \sqrt{\frac{1}{N}
    \sum\limits_{i=1}^N \left(\frac{Y_i^{\text{ref}}-Y_i}{1+\abs{Y_i^{\text{ref}}}} \right)^2},
\end{equation*}
where $Y_i^{\text{ref}}$ and $Y_i$, respectively, denote the reference
solution and the numerical solution at $t_f = 80$ at spatial point $i$
sampled at $N=101$ equally spaced points on the interval $[0,1]$.  The
order of convergence $p$ of the numerical solution was then computed
in the standard way as
\[
p = \frac{\log(\text{error}_1/\text{error}_2 )}{\log(\Delta t_1 / \Delta t_2)},
\]
where the subscripts $1$ and $2$ refer to computations performed using
time steps $\Delta t_1$ and $\Delta t_2$.
\Cref{fig:brusselator_convergence} shows second-order convergence is
observed for all of Strang, \aak1, \aaak2, and \ak52. The number in
the square brackets represents the number of sub-integrations for each
method. From the figure, we see that the alternatives to Strang have
similar (and better) accuracy for a given step size.  However, they
all involve a greater number of sub-integrations and hence are
generally more computationally expensive per step.

\begin{remark}
  The relative sizes of experimental errors observed using Strang,
  \aak1, \aaak2, and \ak52 do not strictly reflect the
  relative sizes of the LEMs presented in \cref{tab:os3_schemes}. This
  is not a surprise because the LEM is only a coarse measure of the
  size of the coefficients that make up the local truncation error. In
  \cite{bernier2020}, the authors observed a smaller error using
  \aaak2 compared to \ak52 when applied to the Vlasov--Maxwell system.

  We further note that by rearranging the order of the operators,
  there are 6 different variants of Strang
  3-splitting~\cref{strang}. All require 5 sub-integrations and have
  similar LEM. However, the cost of solving each sub-integrator is
  different, and the actual error of each of variant is different. We
  consider the following two variants: Strang(2-3-1), whose
  coefficients are given in \cref{tab:strang231} and uses only one
  evaluation of the most time-consuming operator (the diffusion
  operator acting on $T^{[1]}$ and $C^{[1]}$) and Strang(1-3-2), whose
  coefficients are given in \cref{tab:strang132} and generally has the
  smallest actual error for a given step size $\Dt$ for this problem.

  Taking into consideration the number of sub-integrations and the
  size of the error produced, one may hypothesize that only \aak1 is
  competitive versus Strang for this problem, and informal numerical
  experiments with all the methods support this hypothesis.
  Accordingly, we report on the efficiency comparisons only between
  Strang(2-3-1), Strang(1-3-2), and \aak1.
\end{remark}

  \begin{table}[htbp]
  	\centering
  	\begin{tabular}{|c|c|c|c|}
  		\hline
  		$k$&  $\aalpha{k}{1}$ & $\aalpha{k}{2}$ & $\aalpha{k}{3}$ \\
  		\hline
  		1 & $0$ & $0.5$ & $0.5$\\
  		\hline
  		2 & $1$ & $0$ & $0.5$\\
  		\hline
  		3 & $0$ & $0.5$ & $0$\\
  		\hline
  	\end{tabular}
  	\caption{Operator-splitting coefficients of Strang(2-3-1).
  		\label{tab:strang231}}
  \end{table}

  \begin{table}[htbp]
	\centering
	\begin{tabular}{|c|c|c|c|}
		\hline
		$k$&  $\aalpha{k}{1}$ & $\aalpha{k}{2}$ & $\aalpha{k}{3}$ \\
		\hline
		1 & $0.5$ & $0$ & $0.5$\\
		\hline
		2 & $0$ & $1$ & $0.5$\\
		\hline
		3 & $0.5$ & $0$ & $0$\\
		\hline
	\end{tabular}
	\caption{Operator-splitting coefficients of Strang(1-3-2).
		\label{tab:strang132}}
	\end{table}

    To compare the efficiency of Strang(2-3-1), Strang(1-3-2), and
    \aak1, we determine the largest possible step size $\Dt$ can be
    used for each method to achieve MRMS errors of
    $5\%, 4\%, 3\%, 2\%, 1\%$, and $0.5\%$. For each MRMS level, we
    record the minimum wall-clock time needed to perform the
    Strang(2-3-1), Strang(1-3-2), and \aak1 over 10 runs. The
    efficiency tests are performed on a Quad-core Intel Xeon Gold 5122
    CPU 3.60GHz with 48GB of RAM running Ubuntu 18.04 LTS.  We
    observed that Strang(2-3-1) is more efficient than Strang(1-3-2)
    because it evaluates the most expensive sub-integration only once,
    and its accuracy level is comparable to Strang(1-3-2). In
    \cref{tab:efficiency_strang_ak32i}, we present the detailed
    comparison between Strang(2-3-1) and \aak1 at various MRMS
    levels. We notice that \aak1 is on average $9\% - 11\%$ more
    efficient than Strang(2-3-1).

\begin{table}
	\centering
	\begin{tabular}{|c|c|c|c|c|c|}
		\hline
		\multirow{2}{*}{MRMS (\%)} & \multicolumn{2}{c|}{Strang(2-3-1) } & \multicolumn{2}{c|}{\aak1}  & \multirow{2}{*}{Time saved (\%)} \\
		& \multicolumn{1}{c}{$\Dt$} & CPU (s) & \multicolumn{1}{c}{$\Dt$} & CPU (s) & \\
		\hline
		5	& 0.214477212  &	38.4495  & 0.300751880 &	34.7206		& 9.698175529
		\\
		4	& 0.197530864  &	38.7529  & 0.275862069 &	35.1449		& 9.310270973
		\\
		3	& 0.176991150   &	39.6837	 & 0.246913580  & 35.1634	& 11.39082293
		\\
		2	& 0.149532710   &	40.0799  & 0.208333333 &	35.7434		& 10.81963777
		\\
		1	& 0.110041265  &	41.5642  & 0.153256705 &	37.7180	 & 9.253636543
		\\
		0.5	& 0.079522863  & 	44.0782	  & 0.110497238 &	40.0803	& 9.070016471 \\

		\hline
	\end{tabular}
	\caption{Comparison of wall-clock time needed by Strang(2-3-1) and
      \aak1 at various accuracy
      levels. \label{tab:efficiency_strang_ak32i}}
\end{table}

\subsection{Plasma dynamics with kinetic Vlasov--Poisson equations}\label{subsec:VP_results}
Our code is implemented in Fortran 90. Some details of this
implementation can be found in \cite{magdi2009method}. The parameters
used for the ECDI simulations are listed in \cref{tbl:parameters}. The
values of $\alpha_1$, $\alpha_2$, and $\alpha_3$ correspond to the
typical operation regime of Hall thrusters and are also used in
\cite{janhunen2018nonlinear,tavassoli2022nonlinear,janhunen2018evolution}. Also
in this table, ${N_x}$ is the number of cells in the $x$-direction of
the both ion and electron distribution grids; $N_{v_{xe}}$ and
$N_{v_{ze}}$ are the number of cells in the $v_x$- and
$v_z$-directions of the electron grid, respectively; and $N_{v_{xi}}$
is the number of cells in the $v_{x}$-direction of the ion grid. We
note that because the ions are much heavier than the electrons, their
velocity often remains much smaller than that of the electrons during
the simulation.  Because of this, \Cref{eq:ivlasov} is usually solved
on a much coarser grid than \Cref{eq:evlasov}.

\begin{table}[htbp]
	\centering
	\begin{tabular}{|c|c|}
      \hline
      Parameter & Value \\
      \hline
      $\alpha_1$                    & $2.39\times 10^5$    \\
      $\alpha_2$                    & $4.03\times 10^{-4}$ \\
      $\alpha_3$                    & $0.1487$             \\
      ${N_x}$                    & $2048$              \\
      $N_{v_{xe}}$				  & $1200$				 \\
      $N_{v_{ze}}$ 			  & $1200$               \\
      $N_{v_{xi}}$               & $200$              \\
      $L$						& $600$				\\
     \hline
	\end{tabular}

    \caption{Parameters used for the ECDI simulation. \label{tbl:parameters}}
\end{table}

\begin{figure}[htbp]
	\includegraphics[width=\textwidth]{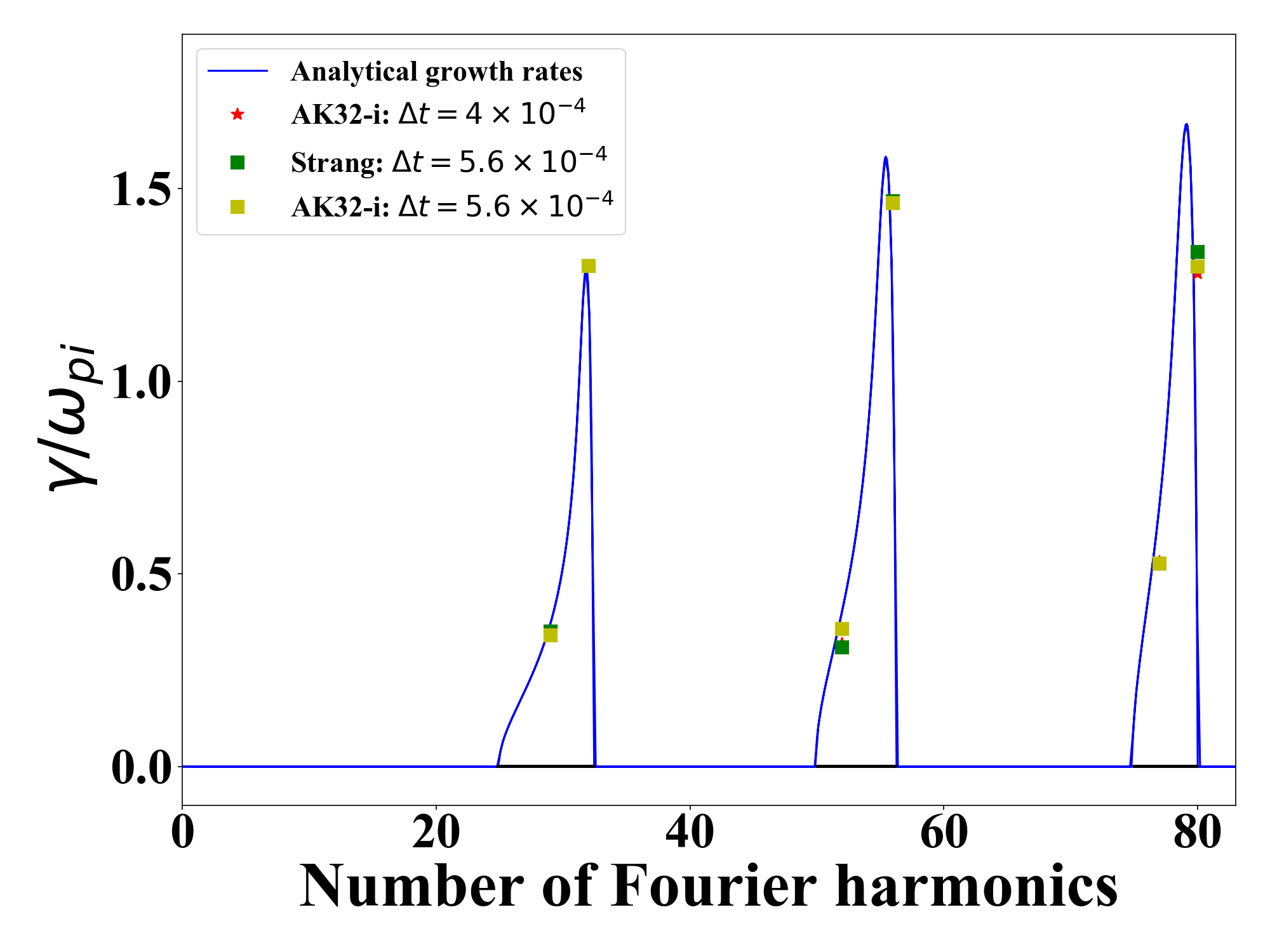}
	\caption{The linear growth rates of the Strang and \aak1
      methods.\label{fig:growth_rates}}
\end{figure}

At the initial stage of the ECDI, called the linear regime, $E_x$
grows exponentially with time. For this growth, different Fourier
harmonics of $E_x$ show different time exponents. These exponents are
referred to as the ``linear growth rates". The linear growth rates can
be calculated analytically by solving a nonlinear algebraic equation
that is derived by applying perturbation methods to the
Vlasov--Poisson system. The details of this derivation can be found in
\cite{gary1993theory}, and a numerical method for solving the
nonlinear algebraic equation is discussed in
\cite{cavalier2013hall}. The linear growth rates can be also estimated
from the simulations, and the results of these estimates can be used
as a measure of the accuracy of the
simulations~\cite{tavassoli2021role,tavassoli2021backward,tavassoli2022nonlinear}. \Cref{fig:growth_rates}
compares the analytical linear growth rates with the growth rates that
are measured numerically by different splitting methods. We can see
that, in all cases, the measured growth rates remain close to each
other and to the analytical growth rates. This means that although both
the Strang and the \aak1 methods provide a reasonable accuracy of the
linear regime, the linear growth rates (and accordingly the short-time
behavior of ECDI) do not provide a suitable measure for distinguishing
the accuracy of these methods.

The linear regime of the instability is followed in time by the
nonlinear regime. In contrast to the linear regime, analysis of the
nonlinear regime is generally performed through simulation.  For
comparing the accuracy of different methods, we use the maximum
relative deviation from energy conservation as a measure of the error
of the simulations. This measure is defined by
\begin{equation*}
  \Delta _t U= 
  \abs{\frac{U(t)-U(0)}{U(0)}}\times 100\%.
\end{equation*}
Because the $\Delta _t U$ is a function of time, its temporal value
may not be a consistent metric for comparing two competing splitting
methods.  Nevertheless, in most applications, a maximum relative
deviation of energy conservation can be tolerated, e.g., 2\%. Here,
although we show the temporal evolution of $\Delta _t U$, we use the
maximum of $\Delta_t U$ over time as a metric to determine solution
accuracy.

\begin{figure}[htbp]
	\includegraphics[width=\textwidth]{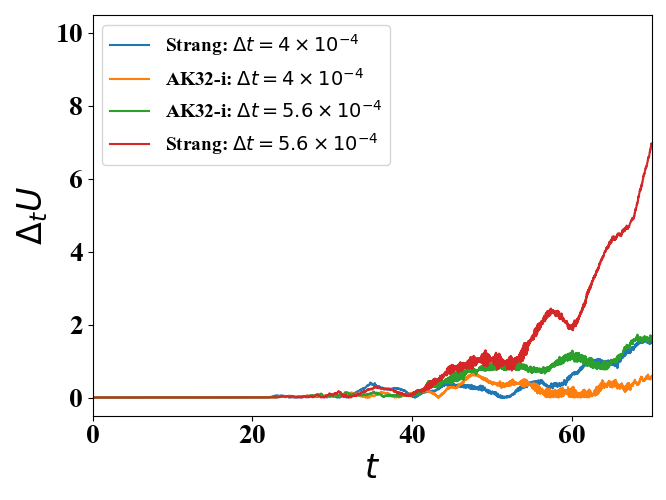}
	\caption{Maximum relative deviation from energy conservation of
      the Strang and \aak1 methods.\label{fig:energy_cons}}
\end{figure}

\Cref{fig:energy_cons} shows $\Delta_t U$ for the Strang and the \aak1
methods. In this figure, the simulation duration is $70$, which
approximately corresponds to the $2\;\mu s$ timespan considered in the
previous studies
\cite{tavassoli2022nonlinear,janhunen2018nonlinear,janhunen2018evolution,lafleur2016theory1}. We
can see that the maximum error is the lowest in \aak1 with
$\Dt=4\times 10^{-4}$. The Strang method with
$\Delta t=4\times 10^{-4}$ has about the same maximum error as the
\aak1 with a time step that is 40\% larger ($\Dt =5.6\times
10^{-4}$). This means that for $\displaystyle \dtfact \approx 1.4$ in
\cref{eq:efficiency}. However, when the time step of the Strang method
is increased by 40\% to $\Delta t=5.6 \times 10^{-4}$, the error jumps
to about 7\%.

For calculating the efficiency gain $\eta$ from \Cref{eq:efficiency},
we also calculated the extra-time fraction $\timeratio$. Compared with
Strang, the \aak1 method needs to evaluate one extra step of
\Cref{eq:vxe_split,eq:vxi_split}. Therefore, one can assume
$\timeratio=\timeratio_{v_{xe}}+\timeratio_{v_{xi}}$, where
$\displaystyle \timeratio_{v_{xe}} =
\frac{\tilde{\tau}^{\aak1}_{\cref{eq:vxe_split}}}{\tstrang}$ and
$\displaystyle \timeratio_{v_{xi}} =\frac{
  \tilde{\tau}^{\aak1}_{\cref{eq:vxi_split}}}{\tstrang}$.  The results
of our measurements of $\timeratio$ are shown in
\cref{table:timeratios} for different sub-integrators of the ECDI. We
note that, because the electron grid is much finer than that of the
ions ($\timeratio_{v_{xe}} \gg \timeratio_{v_{xi}}$), we can neglect
$\timeratio_{v_{xi}}$ arising from \Cref{eq:vxi_split}. Accordingly,
the range of $\timeratio$ is $\timeratio \approx 0.15 - 0.2$, which,
using \Cref{eq:efficiency} and \cref{eq:extracost}, leads to an
efficiency gain of 14\% to 18\% of the \aak1 with respect to
Strang. These estimates are consistent with experimental observations.
\begin{table}[htbp]
	\centering
	\begin{tabular}{|c|c|}
		\hline
		Sub-integrator of
		                   & $\timeratio\times 100$ \\
		\hline
		\Cref{eq:xe_split}  & 5 -- 7            \\
		\Cref{eq:vze_split} & 18 -- 20          \\
		\Cref{eq:vxe_split} & 15 -- 20          \\
		\Cref{eq:xi_split}  & 0.003 -- 0.3      \\
		\Cref{eq:vxi_split} & 0.2 -- 0.6    \\
		\Cref{eq:ef} & 15 -- 25\\
		\hline
	\end{tabular}
	\caption{Estimated ranges for the extra-time fractions
      $\timeratio$ for different sub-integrators of the ECDI
      operators. Measurement performed using 1 to 40 processors and on
      various computers. \label{table:timeratios}}
\end{table}

\section{Discussion and Conclusions}\label{sec:conclusions}

Operator splitting is popular in the numerical solution of
differential equations. The most common framework for splitting is
into two parts; however, there can be advantages to being able to
split a given problem into three parts. The generalization of Strang
splitting to three operators is arguably the most popular 3-splitting
method. Although alternatives to Strang exist, they are less well
known. In this paper, we test some alternative 3-splitting methods on
two applications, the Brusselator and ECDI problems. We demonstrate
worthwhile improvements (10\%--20\%) in efficiency over traditional
Strang splitting, potentially shaving days of computing time off of
real problems that require weeks to simulate (or weeks off of problems
that require months).

In the experiments of the Brusselator problem, we examined various
second-order accurate 3-splitting methods and concluded that \aak1 is
the most efficient method based on overall computation time. Our
experiments have taken into account the LEM error, the number
sub-integrations, and the cost of each sub-integration. We observe
that although all 6 variants of the Strang method exhibit similar
accuracy and require the same number of sub-integrations,
Strang(2-3-1) is the most efficient among the 6 variants because the
accuracy is among the best and the most expensive sub-integration
(that of operator $[1]$ from~\cref{brusselator_1}) is evaluated only
once. We further note that a splitting method may have an optimum LEM
but may require many splitting stages and sub-integrations to solve
one time step. These additional costs may ultimately offset the
efficiency gained from the minimized LEM. Therefore, when choosing an
efficient operator-splitting method, we seek to balance the LEM and
the number of sub-integrations required. This balanced approach is
essential not only to 3-splitting problems but also to general
N-splitting problems. Furthermore, we echo the result in
\cite{auzinger2017} that for some test problems, such as the
Gray–-Scott equations or the Brusselator equations, we may not
experience efficiency gain by splitting the differential equations
into three operators as opposed to two operators. The potential
advantage is that by splitting into three operators, some of the
sub-integrators admit a closed-form solution, which can potentially be
efficiently evaluated.

In the full Vlasov--Poisson equation in 6-dimensional phase space, one
can apply splitting to each direction, and therefore, six operators
need to be solved \cite{kormann2019massively}. Due to their extreme
computational cost, such simulations are generally not feasible at
present. Here, we applied a 3-splitting to a simplified version of the
Vlasov--Poisson system in one spatial direction and two velocity
directions. The 3-splitting helped us to calculate the base of the
characteristic equations explicitly. We also tried a 2-splitting
version of the code, in which the two velocity directions are treated
as only one (rotation) operator. This rotation operator is then
integrated using the method of characteristics, which requires a
2-dimensional B-spline cubic interpolation \cite{magdi2009method}.
For the ECDI problem, we observed that the overall performance
deteriorated compared to the 3-splitting approach. In
\cite{bernier2020}, a new splitting method is proposed for the exact
integration of the rotation operator. Applying the new method to the
Vlasov--Maxwell system, the authors observed a significant increase in
the performance compared to the 2-dimensional spectral
interpolation. Therefore, we hypothesize that this method can also help
increase the performance of our 2-splitting code and propose the
investigation of its impact on the ECDI problem as future work.

In the ECDI problem, we measured $\dtfact\approx 1.4$ and $\timeratio$
in the range $\timeratio\approx 0.15-0.2$, leading to 15\% to 20\%
efficiency gain. This efficiency gain allowed a simulation using 32
processors to complete in about 8 days using \aak1 versus about 10
days using Strang. We have performed other simulations using a
simulation domain four times larger than that used in this work that took
about two months~\cite{tavassoli2023electron,tavassoli2023drift}. In general, 2- and
3-dimensional Vlasov simulations can easily take several months even
using hundreds or thousands of
processors~\cite{tanaka2017multidimensional,kormann2019massively}. Therefore,
the perhaps modest-looking efficiency gain provided by the \aak1
method can lead to several weeks of real time saved for running such
simulations.

To examine the effect of the resolution on $\dtfact$ and $\timeratio$,
we repeated the simulation with the \aak1 method for the $N_{v_{xe}}$
between $800$ to $1200$ and observed that these quantities did not
significantly change. Similarly, running simulations with up to 40 cores did
not have a significant impact on $\timeratio$. We do not expect that
the parallelization or the number of cores to have a significant
impact on $\dtfact$ because this quantity only depends on the accuracy
of a splitting method and not on the number of cores used.

Other than \aak1, we tried other 3-splitting methods, namely \aaak2,
\ak52, and \aaaak11, for the ECDI problem. The general conclusion from
these experiments is that none of these methods are likely to be
competitive with Strang in terms of efficiency. The \aaak2 and \ak52
have four more sub-integrators than Strang, but they approximately
gave the same $\Delta _t U$ as Strang for the same
$\Delta t=4\times 10^{-4}$. Also, the \aaaak11 is a fourth-order
splitting method with 16 sub-integrators more than Strang. Applying
AK11-4 on the ECDI problem and using the values of $\timeratio$ in
\cref{table:timeratios}, we calculate that this method is at least 3
times more expensive per time step than Strang. However, our
experiment with $\Delta t=3\times 4\times 10^{-4}$ showed a much
larger error than Strang with $\Delta t =4\times 10^{-4}$. This may be
explained by the fact that the sub-integrators of our code are based
on cubic spline interpolation, and the error may overshadow the
splitting error of a fourth-order method. Addressing the efficiency of
the splitting methods when higher-order interpolation methods are used
is beyond the scope of the current study.

\backmatter


\bmhead{Acknowledgments}

The authors gratefully acknowledge funding from the Natural Sciences
and Engineering Research Council of Canada under its Discovery Grant
Program (RGPN 2020-04467 (RJS) and RGPN 2022-04482 (AS)) as well as
from the US Air Force Office of Scientific Research FA9550-21-1-0031
(AS).  A. Tavassoli acknowledges the support of Dr. Magdi Shoucri in
developing the semi-Lagrangian code.

\section*{Statements and Declarations}

The authors have no conflict of interest nor competing interests to
declare that are relevant to the content of this article.

\bibliographystyle{bst/sn-mathphys.bst}
\bibliography{sn-bibliography.bib}


\begin{thebibliography}{42}
\ifx \bisbn   \undefined \def \bisbn  #1{ISBN #1}\fi
\ifx \binits  \undefined \def \binits#1{#1}\fi
\ifx \bauthor  \undefined \def \bauthor#1{#1}\fi
\ifx \batitle  \undefined \def \batitle#1{#1}\fi
\ifx \bjtitle  \undefined \def \bjtitle#1{#1}\fi
\ifx \bvolume  \undefined \def \bvolume#1{\textbf{#1}}\fi
\ifx \byear  \undefined \def \byear#1{#1}\fi
\ifx \bissue  \undefined \def \bissue#1{#1}\fi
\ifx \bfpage  \undefined \def \bfpage#1{#1}\fi
\ifx \blpage  \undefined \def \blpage #1{#1}\fi
\ifx \burl  \undefined \def \burl#1{\textsf{#1}}\fi
\ifx \doiurl  \undefined \def \doiurl#1{\url{https://doi.org/#1}}\fi
\ifx \betal  \undefined \def \betal{\textit{et al.}}\fi
\ifx \binstitute  \undefined \def \binstitute#1{#1}\fi
\ifx \binstitutionaled  \undefined \def \binstitutionaled#1{#1}\fi
\ifx \bctitle  \undefined \def \bctitle#1{#1}\fi
\ifx \beditor  \undefined \def \beditor#1{#1}\fi
\ifx \bpublisher  \undefined \def \bpublisher#1{#1}\fi
\ifx \bbtitle  \undefined \def \bbtitle#1{#1}\fi
\ifx \bedition  \undefined \def \bedition#1{#1}\fi
\ifx \bseriesno  \undefined \def \bseriesno#1{#1}\fi
\ifx \blocation  \undefined \def \blocation#1{#1}\fi
\ifx \bsertitle  \undefined \def \bsertitle#1{#1}\fi
\ifx \bsnm \undefined \def \bsnm#1{#1}\fi
\ifx \bsuffix \undefined \def \bsuffix#1{#1}\fi
\ifx \bparticle \undefined \def \bparticle#1{#1}\fi
\ifx \barticle \undefined \def \barticle#1{#1}\fi
\bibcommenthead
\ifx \bconfdate \undefined \def \bconfdate #1{#1}\fi
\ifx \botherref \undefined \def \botherref #1{#1}\fi
\ifx \url \undefined \def \url#1{\textsf{#1}}\fi
\ifx \bchapter \undefined \def \bchapter#1{#1}\fi
\ifx \bbook \undefined \def \bbook#1{#1}\fi
\ifx \bcomment \undefined \def \bcomment#1{#1}\fi
\ifx \oauthor \undefined \def \oauthor#1{#1}\fi
\ifx \citeauthoryear \undefined \def \citeauthoryear#1{#1}\fi
\ifx \endbibitem  \undefined \def \endbibitem {}\fi
\ifx \bconflocation  \undefined \def \bconflocation#1{#1}\fi
\ifx \arxivurl  \undefined \def \arxivurl#1{\textsf{#1}}\fi
\csname PreBibitemsHook\endcsname

\bibitem{ars1997}
\begin{barticle}
\bauthor{\bsnm{Ascher}, \binits{U.M.}},
\bauthor{\bsnm{Ruuth}, \binits{S.J.}},
\bauthor{\bsnm{Spiteri}, \binits{R.J.}}:
\batitle{Implicit-explicit {R}unge--{K}utta Methods for Time-dependent Partial
  Differential Equations}.
\bjtitle{Appl. Numer. Math.}
\bvolume{25}(\bissue{2-3}),
\bfpage{151}--\blpage{167}
(\byear{1997})
\end{barticle}
\endbibitem

\bibitem{Lefever1971}
\begin{barticle}
\bauthor{\bsnm{Lefever}, \binits{R.}},
\bauthor{\bsnm{Nicolis}, \binits{G.}}:
\batitle{Chemical instabilities and sustained oscillations}.
\bjtitle{Journal of Theoretical Biology}
\bvolume{30}(\bissue{2}),
\bfpage{267}--\blpage{284}
(\byear{1971}).
\doiurl{10.1016/0022-5193(71)90054-3}
\end{barticle}
\endbibitem

\bibitem{Strang1968}
\begin{barticle}
\bauthor{\bsnm{Strang}, \binits{G.}}:
\batitle{On the construction and comparison of difference schemes}.
\bjtitle{SIAM Journal on Numerical Analysis}
\bvolume{5}(\bissue{3}),
\bfpage{506}--\blpage{517}
(\byear{1968})
\end{barticle}
\endbibitem

\bibitem{Marchuk1971}
\begin{bchapter}
\bauthor{\bsnm{Marchuk}, \binits{G.I.}}:
\bctitle{On the theory of the splitting-up method}.
In: \beditor{\bsnm{Hubbard}, \binits{B.}} (ed.)
\bbtitle{Numerical Solution of Partial Differential equations–II},
pp. \bfpage{469}--\blpage{500}.
\bpublisher{Academic Press},
\blocation{London}
(\byear{1971})
\end{bchapter}
\endbibitem

\bibitem{auzinger2016practical}
\begin{barticle}
\bauthor{\bsnm{Auzinger}, \binits{W.}},
\bauthor{\bsnm{Hofst{\"a}tter}, \binits{H.}},
\bauthor{\bsnm{Ketcheson}, \binits{D.}},
\bauthor{\bsnm{Koch}, \binits{O.}}:
\batitle{{Practical splitting methods for the adaptive integration of nonlinear
  evolution equations. Part I: Construction of optimized schemes and pairs of
  schemes}}.
\bjtitle{BIT Numerical Mathematics}
\bvolume{57}(\bissue{1}),
\bfpage{55}--\blpage{74}
(\byear{2017})
\end{barticle}
\endbibitem

\bibitem{os_coeff_web}
\begin{botherref}
\oauthor{\bsnm{Auzinger}, \binits{W.}}:
Coefficients of various splitting methods.
\url{http://www.asc.tuwien.ac.at/~winfried/splitting/}
\end{botherref}
\endbibitem

\bibitem{auzinger2017}
\begin{barticle}
\bauthor{\bsnm{Auzinger}, \binits{W.}},
\bauthor{\bsnm{Koch}, \binits{O.}},
\bauthor{\bsnm{Quell}, \binits{M.}}:
\batitle{Adaptive high-order splitting methods for systems of nonlinear
  evolution equations with periodic boundary conditions}.
\bjtitle{Numer. Algorithms}
\bvolume{75}(\bissue{1}),
\bfpage{261}--\blpage{283}
(\byear{2017}).
\doiurl{10.1007/s11075-016-0206-8}
\end{barticle}
\endbibitem

\bibitem{crouseilles2015hamiltonian}
\begin{barticle}
\bauthor{\bsnm{Crouseilles}, \binits{N.}},
\bauthor{\bsnm{Einkemmer}, \binits{L.}},
\bauthor{\bsnm{Faou}, \binits{E.}}:
\batitle{Hamiltonian splitting for the Vlasov--Maxwell equations}.
\bjtitle{Journal of Computational Physics}
\bvolume{283},
\bfpage{224}--\blpage{240}
(\byear{2015})
\end{barticle}
\endbibitem

\bibitem{bernier2020}
\begin{barticle}
\bauthor{\bsnm{Bernier}, \binits{J.}},
\bauthor{\bsnm{Casas}, \binits{F.}},
\bauthor{\bsnm{Crouseilles}, \binits{N.}}:
\batitle{Splitting methods for rotations: application to {V}lasov equations}.
\bjtitle{SIAM J. Sci. Comput.}
\bvolume{42}(\bissue{2}),
\bfpage{666}--\blpage{697}
(\byear{2020}).
\doiurl{10.1137/19M1273918}
\end{barticle}
\endbibitem

\bibitem{hairer2006}
\begin{bbook}
\bauthor{\bsnm{Hairer}, \binits{E.}},
\bauthor{\bsnm{Wanner}, \binits{G.}},
\bauthor{\bsnm{Lubich}, \binits{C.}}:
\bbtitle{Geometric numerical integration: structure-preserving algorithms for
  ordinary differential equations}
vol. \bseriesno{31}.
\bpublisher{Springer},
\blocation{Heidelberg}
(\byear{2006})
\end{bbook}
\endbibitem

\bibitem{Auzinger2014}
\begin{barticle}
\bauthor{\bsnm{Auzinger}, \binits{W.}},
\bauthor{\bsnm{Herfort}, \binits{W.}}:
\batitle{{Local error structures and order conditions in terms of Lie elements
  for exponential splitting schemes}}.
\bjtitle{Opuscula Mathematica}
\bvolume{34}(\bissue{2}),
\bfpage{243}--\blpage{255}
(\byear{2014})
\end{barticle}
\endbibitem

\bibitem{casas2020composition}
\begin{barticle}
\bauthor{\bsnm{Casas}, \binits{F.}},
\bauthor{\bsnm{Escorihuela-Tom{\`a}s}, \binits{A.}}:
\batitle{Composition methods for dynamical systems separable into three parts}.
\bjtitle{Mathematics}
\bvolume{8}(\bissue{4}),
\bfpage{533}
(\byear{2020})
\end{barticle}
\endbibitem

\bibitem{ropp2005}
\begin{barticle}
\bauthor{\bsnm{Ropp}, \binits{D.L.}},
\bauthor{\bsnm{Shadid}, \binits{J.N.}}:
\batitle{Stability of operator splitting methods for systems with indefinite
  operators: reaction-diffusion systems}.
\bjtitle{Journal of Computational Physics}
\bvolume{203}(\bissue{2}),
\bfpage{449}--\blpage{466}
(\byear{2005})
\end{barticle}
\endbibitem

\bibitem{boeuf2018b}
\begin{barticle}
\bauthor{\bsnm{Boeuf}, \binits{J.-P.}},
\bauthor{\bsnm{Garrigues}, \binits{L.}}:
\batitle{{E $\times$ B} electron drift instability in {Hall} thrusters:
  Particle-in-cell simulations vs. theory}.
\bjtitle{Physics of Plasmas}
\bvolume{25}(\bissue{6}),
\bfpage{061204}
(\byear{2018})
\end{barticle}
\endbibitem

\bibitem{charoy2021interaction}
\begin{barticle}
\bauthor{\bsnm{Charoy}, \binits{T.}},
\bauthor{\bsnm{Lafleur}, \binits{T.}},
\bauthor{\bsnm{Laguna}, \binits{A.A.}},
\bauthor{\bsnm{Bourdon}, \binits{A.}},
\bauthor{\bsnm{Chabert}, \binits{P.}}:
\batitle{The interaction between ion transit-time and electron drift
  instabilities and their effect on anomalous electron transport in {Hall}
  thrusters}.
\bjtitle{Plasma Sources Science and Technology}
\bvolume{30}(\bissue{6}),
\bfpage{065017}
(\byear{2021})
\end{barticle}
\endbibitem

\bibitem{sengupta2020mode}
\begin{barticle}
\bauthor{\bsnm{Sengupta}, \binits{M.}},
\bauthor{\bsnm{Smolyakov}, \binits{A.}}:
\batitle{Mode transitions in nonlinear evolution of the electron drift
  instability in a {2D} annular {E $\times$ B} system}.
\bjtitle{Physics of Plasmas}
\bvolume{27}(\bissue{2}),
\bfpage{022309}
(\byear{2020})
\end{barticle}
\endbibitem

\bibitem{asadi2019numerical}
\begin{barticle}
\bauthor{\bsnm{Asadi}, \binits{Z.}},
\bauthor{\bsnm{Taccogna}, \binits{F.}},
\bauthor{\bsnm{Sharifian}, \binits{M.}}:
\batitle{Numerical study of electron cyclotron drift instability: Application
  to Hall thruster}.
\bjtitle{Frontiers in Physics}
\bvolume{7},
\bfpage{140}
(\byear{2019})
\end{barticle}
\endbibitem

\bibitem{hara2020cross}
\begin{barticle}
\bauthor{\bsnm{Hara}, \binits{K.}},
\bauthor{\bsnm{Tsikata}, \binits{S.}}:
\batitle{Cross-field electron diffusion due to the coupling of drift-driven
  microinstabilities}.
\bjtitle{Physical Review E}
\bvolume{102}(\bissue{2}),
\bfpage{023202}
(\byear{2020})
\end{barticle}
\endbibitem

\bibitem{mandal2020cross}
\begin{barticle}
\bauthor{\bsnm{Mandal}, \binits{D.}},
\bauthor{\bsnm{Elskens}, \binits{Y.}},
\bauthor{\bsnm{Lemoine}, \binits{N.}},
\bauthor{\bsnm{Doveil}, \binits{F.}}:
\batitle{Cross-field chaotic transport of electrons by E$\times$ B electron
  drift instability in Hall thruster}.
\bjtitle{Physics of Plasmas}
\bvolume{27}(\bissue{3}),
\bfpage{032301}
(\byear{2020})
\end{barticle}
\endbibitem

\bibitem{janhunen2018nonlinear}
\begin{barticle}
\bauthor{\bsnm{Janhunen}, \binits{S.}},
\bauthor{\bsnm{Smolyakov}, \binits{A.}},
\bauthor{\bsnm{Chapurin}, \binits{O.}},
\bauthor{\bsnm{Sydorenko}, \binits{D.}},
\bauthor{\bsnm{Kaganovich}, \binits{I.}},
\bauthor{\bsnm{Raitses}, \binits{Y.}}:
\batitle{Nonlinear structures and anomalous transport in partially magnetized
  {E$\times$B} plasmas}.
\bjtitle{Physics of Plasmas}
\bvolume{25}(\bissue{1}),
\bfpage{011608}
(\byear{2018})
\end{barticle}
\endbibitem

\bibitem{janhunen2018evolution}
\begin{barticle}
\bauthor{\bsnm{Janhunen}, \binits{S.}},
\bauthor{\bsnm{Smolyakov}, \binits{A.}},
\bauthor{\bsnm{Sydorenko}, \binits{D.}},
\bauthor{\bsnm{Jimenez}, \binits{M.}},
\bauthor{\bsnm{Kaganovich}, \binits{I.}},
\bauthor{\bsnm{Raitses}, \binits{Y.}}:
\batitle{Evolution of the electron cyclotron drift instability in
  two-dimensions}.
\bjtitle{Physics of Plasmas}
\bvolume{25}(\bissue{8}),
\bfpage{082308}
(\byear{2018})
\end{barticle}
\endbibitem

\bibitem{tavassoli2022nonlinear}
\begin{barticle}
\bauthor{\bsnm{Tavassoli}, \binits{A.}},
\bauthor{\bsnm{Smolyakov}, \binits{A.}},
\bauthor{\bsnm{Shoucri}, \binits{M.}},
\bauthor{\bsnm{Spiteri}, \binits{R.J.}}:
\batitle{Nonlinear regimes of the electron cyclotron drift instability in
  {Vlasov} simulations}.
\bjtitle{Physics of Plasmas}
\bvolume{29}(\bissue{3}),
\bfpage{030701}
(\byear{2022})
\end{barticle}
\endbibitem

\bibitem{cheng1976integration}
\begin{barticle}
\bauthor{\bsnm{Cheng}, \binits{C.-Z.}},
\bauthor{\bsnm{Knorr}, \binits{G.}}:
\batitle{The integration of the Vlasov equation in configuration space}.
\bjtitle{Journal of Computational Physics}
\bvolume{22}(\bissue{3}),
\bfpage{330}--\blpage{351}
(\byear{1976})
\end{barticle}
\endbibitem

\bibitem{cheng1977integration}
\begin{barticle}
\bauthor{\bsnm{Cheng}, \binits{C.}}:
\batitle{The integration of the Vlasov equation for a magnetized plasma}.
\bjtitle{Journal of Computational Physics}
\bvolume{24}(\bissue{4}),
\bfpage{348}--\blpage{360}
(\byear{1977})
\end{barticle}
\endbibitem

\bibitem{sonnendrucker1999semi}
\begin{barticle}
\bauthor{\bsnm{Sonnendr{\"u}cker}, \binits{E.}},
\bauthor{\bsnm{Roche}, \binits{J.}},
\bauthor{\bsnm{Bertrand}, \binits{P.}},
\bauthor{\bsnm{Ghizzo}, \binits{A.}}:
\batitle{The semi-Lagrangian method for the numerical resolution of the Vlasov
  equation}.
\bjtitle{Journal of Computational Physics}
\bvolume{149}(\bissue{2}),
\bfpage{201}--\blpage{220}
(\byear{1999})
\end{barticle}
\endbibitem

\bibitem{crouseilles2010conservative}
\begin{barticle}
\bauthor{\bsnm{Crouseilles}, \binits{N.}},
\bauthor{\bsnm{Mehrenberger}, \binits{M.}},
\bauthor{\bsnm{Sonnendr{\"u}cker}, \binits{E.}}:
\batitle{Conservative semi-Lagrangian schemes for Vlasov equations}.
\bjtitle{Journal of Computational Physics}
\bvolume{229}(\bissue{6}),
\bfpage{1927}--\blpage{1953}
(\byear{2010})
\end{barticle}
\endbibitem

\bibitem{crouseilles2009forward}
\begin{barticle}
\bauthor{\bsnm{Crouseilles}, \binits{N.}},
\bauthor{\bsnm{Respaud}, \binits{T.}},
\bauthor{\bsnm{Sonnendr{\"u}cker}, \binits{E.}}:
\batitle{A forward semi-Lagrangian method for the numerical solution of the
  Vlasov equation}.
\bjtitle{Computer Physics Communications}
\bvolume{180}(\bissue{10}),
\bfpage{1730}--\blpage{1745}
(\byear{2009})
\end{barticle}
\endbibitem

\bibitem{qiu2010conservative}
\begin{barticle}
\bauthor{\bsnm{Qiu}, \binits{J.-M.}},
\bauthor{\bsnm{Christlieb}, \binits{A.}}:
\batitle{A conservative high order semi-Lagrangian WENO method for the Vlasov
  equation}.
\bjtitle{Journal of Computational Physics}
\bvolume{229}(\bissue{4}),
\bfpage{1130}--\blpage{1149}
(\byear{2010})
\end{barticle}
\endbibitem

\bibitem{besse2008convergence}
\begin{barticle}
\bauthor{\bsnm{Besse}, \binits{N.}},
\bauthor{\bsnm{Mehrenberger}, \binits{M.}}:
\batitle{Convergence of classes of high-order semi-Lagrangian schemes for the
  Vlasov--Poisson system}.
\bjtitle{Mathematics of computation}
\bvolume{77}(\bissue{261}),
\bfpage{93}--\blpage{123}
(\byear{2008})
\end{barticle}
\endbibitem

\bibitem{ghizzo2003non}
\begin{barticle}
\bauthor{\bsnm{Ghizzo}, \binits{A.}},
\bauthor{\bsnm{Huot}, \binits{F.}},
\bauthor{\bsnm{Bertrand}, \binits{P.}}:
\batitle{A non-periodic 2D semi-Lagrangian Vlasov code for laser--plasma
  interaction on parallel computer}.
\bjtitle{Journal of Computational Physics}
\bvolume{186}(\bissue{1}),
\bfpage{47}--\blpage{69}
(\byear{2003})
\end{barticle}
\endbibitem

\bibitem{coulaud1999parallelization}
\begin{barticle}
\bauthor{\bsnm{Coulaud}, \binits{O.}},
\bauthor{\bsnm{Sonnendr{\"u}cker}, \binits{E.}},
\bauthor{\bsnm{Dillon}, \binits{E.}},
\bauthor{\bsnm{Bertrand}, \binits{P.}},
\bauthor{\bsnm{Ghizzo}, \binits{A.}}:
\batitle{Parallelization of semi-Lagrangian Vlasov codes}.
\bjtitle{Journal of Plasma Physics}
\bvolume{61}(\bissue{3}),
\bfpage{435}--\blpage{448}
(\byear{1999})
\end{barticle}
\endbibitem

\bibitem{magdi2009method}
\begin{bchapter}
\bauthor{\bsnm{Shoucri}, \binits{M.}}:
\bctitle{The method of characteristics for the numerical solution of hyperbolic
  differential equations}.
In: \beditor{\bsnm{Baswell}, \binits{A.R.}} (ed.)
\bbtitle{Advances in Mathematics Research, Volume 8},
pp. \bfpage{1}--\blpage{87}.
\bpublisher{Nova Science Publishers, Inc.},
\blocation{New York}
(\byear{2009}).
\bcomment{Chap. 1}
\end{bchapter}
\endbibitem

\bibitem{staniforth1991semi}
\begin{barticle}
\bauthor{\bsnm{Staniforth}, \binits{A.}},
\bauthor{\bsnm{C{\^o}t{\'e}}, \binits{J.}}:
\batitle{{Semi-Lagrangian} integration schemes for atmospheric models—A
  review}.
\bjtitle{Monthly weather review}
\bvolume{119}(\bissue{9}),
\bfpage{2206}--\blpage{2223}
(\byear{1991})
\end{barticle}
\endbibitem

\bibitem{gary1993theory}
\begin{bbook}
\bauthor{\bsnm{Gary}, \binits{S.P.}}:
\bbtitle{Theory of space plasma microinstabilities}
vol. \bseriesno{7}.
\bpublisher{Cambridge University Press},
\blocation{Cambridge}
(\byear{1993})
\end{bbook}
\endbibitem

\bibitem{cavalier2013hall}
\begin{barticle}
\bauthor{\bsnm{Cavalier}, \binits{J.}},
\bauthor{\bsnm{Lemoine}, \binits{N.}},
\bauthor{\bsnm{Bonhomme}, \binits{G.}},
\bauthor{\bsnm{Tsikata}, \binits{S.}},
\bauthor{\bsnm{Honore}, \binits{C.}},
\bauthor{\bsnm{Gresillon}, \binits{D.}}:
\batitle{Hall thruster plasma fluctuations identified as the {E $\times$B}
  electron drift instability: Modeling and fitting on experimental data}.
\bjtitle{Physics of Plasmas}
\bvolume{20}(\bissue{8}),
\bfpage{082107}
(\byear{2013})
\end{barticle}
\endbibitem

\bibitem{tavassoli2021role}
\begin{barticle}
\bauthor{\bsnm{Tavassoli}, \binits{A.}},
\bauthor{\bsnm{Chapurin}, \binits{O.}},
\bauthor{\bsnm{Jimenez}, \binits{M.}},
\bauthor{\bsnm{Papahn~Zadeh}, \binits{M.}},
\bauthor{\bsnm{Zintel}, \binits{T.}},
\bauthor{\bsnm{Sengupta}, \binits{M.}},
\bauthor{\bsnm{Cou{\"e}del}, \binits{L.}},
\bauthor{\bsnm{Spiteri}, \binits{R.J.}},
\bauthor{\bsnm{Shoucri}, \binits{M.}},
\bauthor{\bsnm{Smolyakov}, \binits{A.}}:
\batitle{The role of noise in {PIC} and {Vlasov} simulations of the {Buneman}
  instability}.
\bjtitle{Physics of Plasmas}
\bvolume{28}(\bissue{12}),
\bfpage{122105}
(\byear{2021})
\end{barticle}
\endbibitem

\bibitem{tavassoli2021backward}
\begin{barticle}
\bauthor{\bsnm{Tavassoli}, \binits{A.}},
\bauthor{\bsnm{Shoucri}, \binits{M.}},
\bauthor{\bsnm{Smolyakov}, \binits{A.}},
\bauthor{\bsnm{Papahn~Zadeh}, \binits{M.}},
\bauthor{\bsnm{Spiteri}, \binits{R.J.}}:
\batitle{Backward waves in the nonlinear regime of the Buneman instability}.
\bjtitle{Physics of Plasmas}
\bvolume{28}(\bissue{2}),
\bfpage{022307}
(\byear{2021})
\end{barticle}
\endbibitem

\bibitem{lafleur2016theory1}
\begin{barticle}
\bauthor{\bsnm{Lafleur}, \binits{T.}},
\bauthor{\bsnm{Baalrud}, \binits{S.}},
\bauthor{\bsnm{Chabert}, \binits{P.}}:
\batitle{Theory for the anomalous electron transport in Hall effect thrusters.
  I. Insights from particle-in-cell simulations}.
\bjtitle{Physics of Plasmas}
\bvolume{23}(\bissue{5}),
\bfpage{053502}
(\byear{2016})
\end{barticle}
\endbibitem

\bibitem{kormann2019massively}
\begin{barticle}
\bauthor{\bsnm{Kormann}, \binits{K.}},
\bauthor{\bsnm{Reuter}, \binits{K.}},
\bauthor{\bsnm{Rampp}, \binits{M.}}:
\batitle{A massively parallel semi-Lagrangian solver for the six-dimensional
  Vlasov--Poisson equation}.
\bjtitle{The International Journal of High Performance Computing Applications}
\bvolume{33}(\bissue{5}),
\bfpage{924}--\blpage{947}
(\byear{2019})
\end{barticle}
\endbibitem

\bibitem{tavassoli2023electron}
\begin{botherref}
\oauthor{\bsnm{Tavassoli}, \binits{A.}},
\oauthor{\bsnm{Papahn~Zadeh}, \binits{M.}},
\oauthor{\bsnm{Smolyakov}, \binits{A.}},
\oauthor{\bsnm{Shoucri}, \binits{M.}},
\oauthor{\bsnm{Spiteri}, \binits{R.J.}}:
The electron cyclotron drift instability: A comparison of particle-in-cell and
  continuum Vlasov simulations.
Physics of Plasmas
\textbf{30}(3)
(2023)
\end{botherref}
\endbibitem

\bibitem{tavassoli2023drift}
\begin{botherref}
\oauthor{\bsnm{Tavassoli}, \binits{A.}}:
Drift instabilities, anomalous transport, and heating in low-temperature
  plasmas.
PhD thesis,
University of Saskatchewan
(2023)
\end{botherref}
\endbibitem

\bibitem{tanaka2017multidimensional}
\begin{barticle}
\bauthor{\bsnm{Tanaka}, \binits{S.}},
\bauthor{\bsnm{Yoshikawa}, \binits{K.}},
\bauthor{\bsnm{Minoshima}, \binits{T.}},
\bauthor{\bsnm{Yoshida}, \binits{N.}}:
\batitle{Multidimensional Vlasov--Poisson simulations with high-order
  monotonicity- and positivity-preserving schemes}.
\bjtitle{The Astrophysical Journal}
\bvolume{849}(\bissue{2}),
\bfpage{76}
(\byear{2017})
\end{barticle}
\endbibitem

\end{thebibliography}


\end{document}